\documentclass[Journal]{IEEEtran}
\usepackage{amsmath,amssymb,bm,amsthm}
\usepackage{algorithm,algorithmic,ulem}
\usepackage{subfig}
\usepackage{graphicx}
\usepackage{booktabs}
\usepackage[english]{babel}

\usepackage{comment,physics}
\usepackage{stackengine}
\usepackage{lineno}
\usepackage{xcolor}
\usepackage{array,multirow}
\usepackage{makecell}
\usepackage{cite}

\usepackage{scalerel}

\def\dashint{\,\ThisStyle{\ensurestackMath{%
\stackinset{c}{.2\LMpt}{c}{.5\LMpt}{\SavedStyle-}{\SavedStyle\phantom{\int}}}%
        \setbox0=\hbox{$\SavedStyle\int\,$}\kern-\wd0}\int}

\newcommand{\reviewerone}[1]{\textcolor{black!50!black}{#1}}

\begin{document}
\bstctlcite{IEEEexample:BSTcontrol}
    \title{Iso-geometric Integral Equation Solvers and their Compression via Manifold Harmonics}

  \author{A. M. A. Alsnayyan,~\IEEEmembership{Student Member,~IEEE,}
      and~B.~Shanker,~\IEEEmembership{Fellow,~IEEE}

  \thanks{The authors acknowledge computing support from the HPC Center at Michigan State University, financial support from NSF via CMMI-1725278 and US Air Force Research Laboratory under contracts FA8650-19-F-1747 and FA8650-20-C-1132.}
  
  \thanks{The  authors  are  with  the  Department  of  Electrical  and  Computer Engineering, Michigan State University, East Lansing, MI 48824-1226 USA. (e-mail: alsnayy1@msu.edu).}
}

\markboth{IEEE TRANSACTIONS ON ANTENNAS AND PROPAGATION
}{Alsnayyan \MakeLowercase{\textit{et al.}}: Well conditioned Integral Equations}

\maketitle

\begin{abstract}
The state of art of electromagnetic integral equations has seen significant growth over the past few decades, overcoming some of the fundamental bottlenecks: computational complexity, low frequency and dense discretization breakdown, preconditioning, and so on. Likewise, the community has seen extensive investment in development of methods for higher order analysis, in both geometry and physics. Unfortunately, these standard geometric descriptors are \reviewerone{continuous, but their normals are discontinuous} at the boundary between \reviewerone{triangular tessellations  of control nodes, or patches,} with a few exceptions; as a result, one needs to define additional mathematical infrastructure to define physical basis sets for vector problems. In stark contrast, the geometric representation used for design \reviewerone{are second order differentiable almost} everywhere on the surfaces. Using these description for analysis opens the door to several possibilities, and is the area we explore in this paper. Our focus is on Loop subdivision based isogeometric methods. In this paper, our goals are two fold: (i) development of computational infrastructure for isogeometric analysis of electrically large simply connected objects, and (ii) to introduce the notion of manifold harmonics transforms and its utility in computational electromagnetics. Several results highlighting the efficacy of these two methods are presented.
\end{abstract}

\begin{IEEEkeywords}
Integral equations, subdivision surfaces, iso-geometric methods, higher order, manifold harmonics, Fast  multipole  method
\end{IEEEkeywords}

\IEEEpeerreviewmaketitle

\section{Introduction}

\IEEEPARstart{O}{ver} the past six decades, the state of the art boundary integral equation solvers have grown by leaps and bounds to become a powerful tool for electromagnetic analysis. A sequence of advancements have enabled this transition, starting from the development of integral equations (see \cite{Peterson1997b} and references therein for a more complete historical background), to methods to appropriately discretize them \cite{RWG}, to higher order representations \cite{pete_high_order}, to overcoming computational bottlenecks \cite{Wandzuraz1993,Song1997,dault2016mixed,vikram2009novel}, to well conditioned formulations \cite{Turc_IE,Chew2014,Rokhlin_IE}, and more recently, to preconditioning techniques \cite{andriulli_precond,simon_phdthesis}. However, despite the significant recent progress made, the technological drivers demand a more sophisticated and more feature rich solver, \emph{albeit} at reduced cost. 

Computational analysis typically proceeds in three stages; (a) construct a geometric model using a computer aided design (CAD) tool, (b) define a discrete representation of said geometry, and (c) finally, choosing a representation of the physics on the discrete representation of the geometry. Geometry is typically represented using bi-variate splines (Bezier splines, B-splines, or non-uniform rational B-splines (NURBS)) that can provide higher order continuity on the surface. From this surface representation, a mesh is generated that typically provides low order continuity on the manifold. As an example, piecewise flat Lagrangian elements are $C^0$, i.e., continuous at interfaces between \reviewerone{triangular tesselations of control nodes, or patches,} but with discontinuous normals. Furthermore, even higher order meshes are higher order within a patch/subdomain, but still $C^0$ across patches. As a result, basis functions defined on these meshes must impose additional constraints. In this framework, a number of different approaches to electromagnetic analysis tools have been developed, including: RWG basis sets \cite{RWG}, its higher order variants \cite{pete_high_order}, Buffa-Christansen basis \cite{BC_basis}. \reviewerone{In addition}, there exists an in-depth analysis and study into computational bottlenecks such as ill-conditioning, low-frequency breakdown, dense-mesh breakdown, topological breakdown, etc \cite{Rokhlin_IE,Kleinman,andriulli_precond}. 

Two more relatively recent methods take a different approach; they still seek to obtain a higher order parameterization of the geometry and thereby, higher order basis for physics. The first overcomes item (a) above, and directly models the object using higher order polynomials \cite{Notaros}. Another approach, the generalized method of moments (GMM), starts with (b) and builds a framework that accommodates both large ($>4\lambda$) and small patches as well as different functions on each patch \cite{DaultGMM,JieDaultShankerChapterSubd} all stitched together within a partition of unity framework. This is done using a non-watertight \reviewerone{representation} of standard meshes. Other methods rely on different techniques to enrich function spaces to represent physics (for instance, macro-basis sets \cite{MBF_Vion}). All seek to achieve \reviewerone{an} efficient representation of geometry, physics, or both. 

An alternative approach that is gaining currency is equipped with the infrastructure to do physics using the same basis function used to construct the geometry; this is known as \reviewerone{isogeometric} analysis (IGA). The advantages of such an approach are as follows: they (a) eliminate the error in translating between geometry and the mesh; (b) the number of degrees of freedom is limited to that used for geometry representation which is \emph{significantly smaller} than a corresponding mesh; and (c) the rules used for adaptation and refinement are identical for both geometry and physics; a vivid illustration can be found in \cite{Buffa,Simpson,JieDaultShankerChapterSubd,Jie}. One must highlight that in isogeometric methods, basis functions are co-located on control nodes used to describe the geometry. This is in contrast with parametric methods that require additional infrastructure--for an example of using subdivision for geometry and GMM basis sets, see \cite{JieDaultShankerChapterSubd}. 

The genesis of IGA methods started with using NURBS for solid mechanics \cite{Hughes_IGA}, and more recently, in electromagnetics \reviewerone{\cite{Buffa,iga_EM}} and acoustics \cite{Simpson}. \reviewerone{NURBS are geometric descriptions that are topologically either a disk, a tube or a torus. As a result, one of the major difficulties arise with NURBs that the patches have to be seamlessly sewed together in order to handle complex surfaces which is time-consuming and complicated. Furthermore, stitching together these patches can result in surfaces that are not watertight and sometimes discontinuous. These complexities have proven to be quite a hindrance when handling complex geometries \cite{Bazilevs2010}. Other modalities that have gained currency in geometry representation are T-splines and Loop subdivision. While T-splines have been used in an IGA setting (see \cite{Bazilevs2010,Doerfel2010} and references therein) our focus in this paper will be on Loop subdivision. }

Loop subdivision have been extremely popular in the computer graphics industry due to the ease with which one can represent complex topologies, its scalability, inherently multiresolution features, efficiency and ease of implementation. More importantly, the surface representation is $C^2$, \reviewerone{or continuous twice differentiable surface,} almost everywhere making it an attractive candidate for defining physical basis sets as it avoids the requirement of defining additional mathematical framework that is commonplace in other low order basis set \cite{Kleinman,RWG,high_order_RWG}. There has been a concerted effort to develop IGA methods on subdivision surfaces in a number of fields, including electromagnetics \cite{Li2020,Jie,Fu2017GeneralizedDS}, acoustics \cite{abdel_acoustics,Simpson} and shape reconstruction/optimization \cite{alsnayyan2021laplacebeltrami,FEM_multires,shape_FEM_multi,Taka_photo,acoustics_topo_opt_subd}.

This paper builds on our earlier body of work on Loop subdivision based IGA for the electric field integral equations \cite{Jie} and construction of Debye sources \cite{Fu2017GeneralizedDS}. In both these cases, the objects analyzed were simply connected \emph{and} electrically small. \reviewerone{Further, they only discreized the electric field integral equation.} The key bottleneck is the number of quadrature points required to evaluate all necessary inner products on higher order geometry (4$^{th}$) and 3$^{rd}$ order basis. A principal goal of this paper is to alleviate this bottleneck for \emph{all methods that use higher order surface representation and higher order basis for physics. It is illustrated here for subdivision basis}. To do so, we exploit wideband multilevel fast multipole algorithm to evaluate \emph{all} interactions (self, near, and far) with leaf boxes as small as 0.025$\lambda$. Furthermore, we pair this approach with a well conditioned combined field integral equation to analyze  objects  as large as 120$\lambda$.

Next, we introduce manifold harmonic basis (MHB) for field computation. These basis are the eigenfunctions of the Laplace Beltrami Operator (LBO) \cite{Italian_LBO_book} and are computed using finite element on the manifold. MHB is tantamount to Fourier basis on the manifold \cite{Levy_geom_process} \reviewerone{and analogously is equipped with a manifold harmonic transform (MHT)}.
It has found numerous applications, ranging from shape analysis \cite{shape_analy,LBO_DNA}, dimensionality reduction with spectral embeddings \cite{LBO_reduction,Richard_LBO}, 
medical imagining applications\cite{laplace_imag,Moo_LBO}, and shape reconstruction \cite{alsnayyan2021laplacebeltrami}. In this paper, we explore the applicability of MHB for electromagnetic analysis, specifically to compress systems resulting from discretization of boundary integral equations in electromagnetics, and demonstrate its numerous benefits. \reviewerone{\emph{What we do not address}, and is outside the scope of this paper, is the cost of applying these transformation, remediation of cost and the other benefits that arise from this transformation;} these topics will be addressed in subsequent papers and the direction of our research on these issues is alluded to in summary section of this paper.

\section{Problem statement}

We consider the analysis of scattered fields $\{ \textbf{E}^{s},\textbf{H}^{s}\}$, from a perfect electrically conducting (PEC) object
$\Omega$, due to fields $\{ \textbf{E}^{i},\textbf{H}^{i}\}$ incident on its boundary $\Gamma \in \Omega$. It is assumed that this surface is equipped with a unique outward pointing normal denoted by $\hat{\textbf{n}}(\textbf{r})$, $\textbf{r} \in \Gamma$. The region external to this volume $\{\mathbb{R}^{3}\setminus \Omega\}$ is occupied by free space. The scattered field at $\textbf{r} \in \{\mathbb{R}^{3}\setminus\Omega\}$ can be obtained using equivalence theorems leading to the following:

\begin{align}
    \begin{split}
            \hat{\textbf{n}}(\textbf{r}) \times \textbf{E}^{s}(\textbf{r}) &= \mathcal{T}_{\kappa}\circ \textbf{J}(\textbf{r}),\\
            \hat{\textbf{n}}(\textbf{r}) \times \textbf{H}^{s}(\textbf{r}) &= \mathcal{K}_{\kappa}\circ \textbf{J}(\textbf{r}),
    \end{split}
\end{align}
where,
\begin{subequations}
\begin{align}
\begin{split}
         \mathcal{T}_{\kappa} \circ \textbf{J}(\textbf{r}) &=  -j\eta \kappa\hat{\textbf{n}}(\textbf{r}) \times \int_{\Gamma} G_{\kappa}(\textbf{r},\textbf{r}')\cdot \textbf{J}(\textbf{r}')d\textbf{r}' \\  &+ \frac{\eta}{j\kappa}\hat{\textbf{n}}(\textbf{r}) \times \nabla\int_{\Gamma}  G_{\kappa}(\textbf{r},\textbf{r}') \nabla' \cdot \textbf{J}(\textbf{r}')d\textbf{r}',
\end{split}
\end{align}
\reviewerone{
\begin{equation}
    \mathcal{K}_{\kappa} \circ \textbf{J}(\textbf{r}) =  \hat{\textbf{n}}(\textbf{r}) \times \dashint_{\Gamma} \nabla G_{\kappa}(\textbf{r},\textbf{r}')\cdot \textbf{J}(\textbf{r}')d\textbf{r}', 
\end{equation}
}
\end{subequations}
where $G_{\kappa}(\textbf{r},\textbf{r}') = \mbox{exp}[-j\kappa\abs{\textbf{r}-\textbf{r}'}]/(4 \pi \abs{\textbf{r}-\textbf{r}'})$, $\kappa$ is the free space wavenumber, $\eta$ is the free space impedance, \reviewerone{ $\mathcal{K}_{\kappa}$ is taken in the Cauchy principal value sense}, and $\textbf{J}(\textbf{r}')$ is the equivalent current that is induced on the surface. In the above expressions, and what follows, we assume and suppress $\mbox{exp}[j\omega t]$ time dependence. Using the above equations, one may prescribe the requisite electric field and magnetic field integral equations (EFIE/MFIE) as 
\begin{subequations}
\begin{align}
    \begin{split}
            \mbox{\textbf{EFIE: }} &= \hat{\textbf{n}}(\textbf{r}) \times \hat{\textbf{n}}(\textbf{r})  \times \left (\textbf{E}^{i}(\textbf{r}) + \textbf{E}^{s}(\textbf{r}) \right ) = 0, 
    \end{split}\\
    \begin{split}
    \mbox{ \textbf{ MFIE: }} &= \hat{\textbf{n}}(\textbf{r})  \times \left (\textbf{H}^{i}(\textbf{r}) + \textbf{H}^{s}(\textbf{r}) \right) = \frac{\textbf{J}(\textbf{r})}{2} .
    \end{split}
\end{align}
\end{subequations}
Independently, these equations suffer from non-unique solutions at so-called irregular frequencies, but their linear combination yields a uniquely solvable formulation throughout the frequency spectrum denoted as the combined field integral equation (CFIE):

\begin{equation}
\begin{split}
\label{eq:CFIE}
                (1-\alpha)\left (\reviewerone{\frac{\mathcal{I}}{2}} - \mathcal{K}_{\kappa} \right ) \circ \textbf{J} (\vb{r}) &+ \alpha \hat{\textbf{n}} \times \mathcal{T}_{\kappa}\circ \textbf{J} (\vb{r}) =\\  &(1-\alpha)\hat{\textbf{n}} \times \textbf{H}^{i} - \alpha \hat{\textbf{n}} \times \hat{\textbf{n}} \times \textbf{E}^{i} ,
\end{split} 
\end{equation}
where $\alpha$ is a positive constant. It is well known that these integral equations suffer from several breakdowns (low frequency, dense mesh, topology, etc.) \cite{Turc_IE,Yass_IE,Rokhlin_IE}. There has been an extensive body of literature addressing these bottlenecks \cite{Jie,andriulli_precond}. In particular, in \cite{Jie,JieDaultShankerChapterSubd,Fu2017GeneralizedDS} the following has been demonstrated for the EFIE: for simply connected objects, employing an isogeometric framework, it is then possible to create a basis that completely satisfy Helmholtz decomposition and this basis set can be used in a Calder\'{o}n setting. While this overcomes a number of problems, a regularized CFIE formulations is still necessary to overcome the non-uniqueness problem. In what follows, we detail a regularized CFIE.    

\subsection{Regularized Combined Field Integral Equations (CFIER)}

A regularized reformulation of (\ref{eq:CFIE}) is the \reviewerone{CFIER} written as follows: 
\reviewerone{
\begin{equation}
    \label{eq:CCFIE-R}
    \left (\reviewerone{\frac{\mathcal{I}}{2}} - \mathcal{K}_{\kappa} \right ) \circ \textbf{J} + \mathcal{R}_{\kappa} \circ \mathcal{T}_{\kappa} \circ \mathbf{J} = \hat{\textbf{n}} \cross \textbf{H}^{i} - \mathcal{R}_{\kappa}\circ (\hat{\textbf{n}} \cross \textbf{E}^{i}). 
\end{equation}
}
 Here, $\mathcal{R}_{\kappa}$ is chosen as a regularizing operator for $\mathcal{T}_{\kappa}$ such that the integral operators on the left hand side of (\ref{eq:CCFIE-R}) are second kind Fredholm operators. Typically, the construction of the regularizing operators is based on \reviewerone{Calder\'{o}n} identities and complexification techniques. Operator $\mathcal{R}_{\kappa}$ have been proposed and analyzed in the literature \reviewerone{ \cite{Rokhlin_IE,Turc_IE,Yass_IE,DARBAS2006834,mitharwal2014multiplicative,dely2019preconditioning}}. 

In particular, we choose the regularization operators provided in \cite{Yass_IE}. This formulation was found to showcase the superior performance of solvers based on the novel \reviewerone{Calder\'{o}n}-Complex CFIER (CC-CFIER) formulations that involve the boundary integral operators
\begin{equation}
    \label{eq:CCCFIER}
    \left  (\reviewerone{\frac{\mathcal{I}}{2}} - \mathcal{K}_{\kappa} \right) \circ \textbf{J}  - 2\mathcal{T}_{\kappa'} \circ \mathcal{T}_{\kappa} \circ \mathbf{J} = \hat{\textbf{n}} \cross \textbf{H}^{i} + 2\mathcal{T}_{\kappa'} \circ (\hat{\textbf{n}} \cross \textbf{E}^{i}),
\end{equation}
where \reviewerone{$\kappa' = \kappa - j0.4\varsigma^{2/3}\kappa^{1/3}$} and $\varsigma$ is the maximum of the absolute values of mean curvatures on surface $\Gamma$.

To solve \eqref{eq:CCCFIER} we will (i) represent the surface of the scatterer using \reviewerone{isogeometric} Loop subdivision basis sets, (ii) represent the currents on the surface using the \emph{same} basis set, and (iii) validate solutions to these integral equations solved using this procedure. Next, we discuss these in sequence.

\section{Subdivision Surfaces and Functions} 
 
\begin{figure}[!t]
  \centering
  \includegraphics[width=3cm]{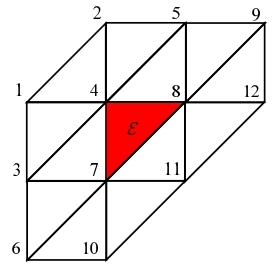}
  \caption{Regular triangular patch defined by its 1-ring vertices.
  \label{fig:1ring}}  
\end{figure}

\begin{figure*}[!t]
\hspace*{-1.2cm}              
\centering
\subfloat[$H_{1}$]{
    \label{fig:bumpy_cube_1}
	\includegraphics[width=0.25\textwidth]{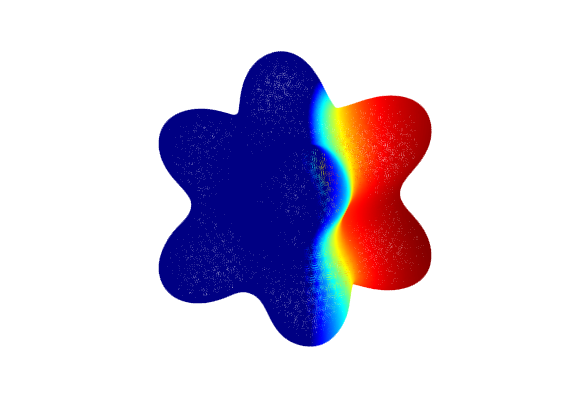} }
	\subfloat[$H_{9}$]{
    \label{fig:bumpy_cube_9}
	\includegraphics[width=0.25\textwidth]{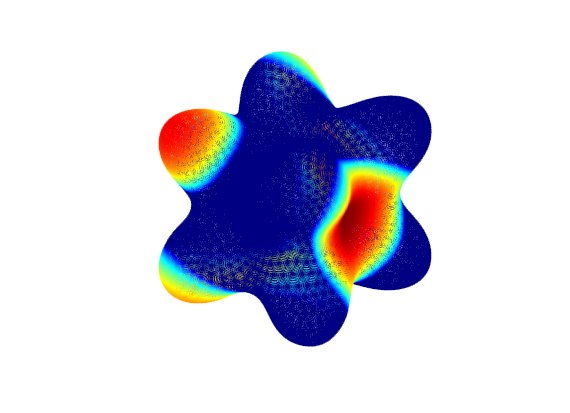} }
\subfloat[$H_{500}$]{
    \label{fig:bumpy_cube_500}
	\includegraphics[width=0.25\textwidth]{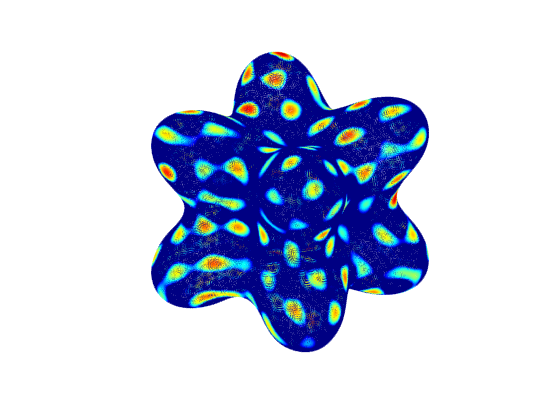} } 
%
\subfloat[$H_{1500}$]{
    \label{fig:bumpy_cube_1500}
	\includegraphics[width=0.25\textwidth]{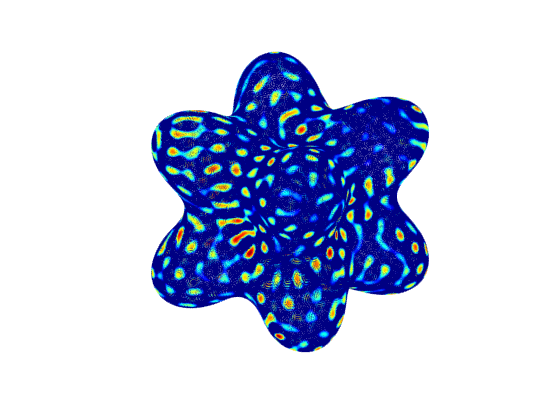} } 
 \\
\caption{ A select few MHs of the bumpy cube. (a) $H_{1}$. (b) $H_{9}$. (c) $H_{500}$. (d) $H_{1500}$.}
\label{fig:bumpy_cube_MHs}
\end{figure*}    

In this Section, we provide a brief overview of Loop subdivision as an isogeometric tool; information provided is purely for completeness and omits details that can be found in \cite{Ciraksubd,stam1998evaluation,loop1987smooth,JieDaultShankerChapterSubd,abdel_acoustics,IGA_JCP,Rumpf_IGA} and references therein. Let $T^{k}$ denote a $k$-th refined control mesh, with vertices $V^{k} := \{\textbf{v}_{i}, i = 1,\ldots, N_{v}\}$ and triangular faces $P^{k} := \{\textbf{p}_{i}, i = 1,\ldots, N_{f}\}$. In short, we can represent a $C^{2}$ (almost everywhere) smooth limit surface $\Gamma$, through an infinite number of iterative refinements of the control mesh $T^{0}$, following the loop subdivision scheme \cite{subdivision}. In practice, this prescription is not followed.  There exists closed form expressions for computing the limit surface $\Gamma$ for a given control mesh $T^{k}$ in terms of quantities defined on the given control mesh \cite{stam1998evaluation}. Assume that a subdivision surface admits a natural parameterization of the surface $\Gamma$ in terms of the barycentric coordinates defined on each face $\epsilon \in P^{k}$, for some {\it k}. We begin by considering any patch $\epsilon \in P^{k}$ for some {\it k}, as depicted in Fig.~\ref{fig:1ring}. We define the 0-ring of a patch (triangle) as the vertices that belong to the patch, and the 1-ring as the set of all vertices, $n_{v}$, that can be reached by traversing no more than two edges, as shown in Fig.~\ref{fig:1ring}. We define the regularity of the triangle by the characterization of its vertices' valence (0-ring); the valence of a given vertex is the number of edges incident on itself. A vertex is considered regular if its valence is equal to 6, otherwise, it is called an irregular or extraordinary vertex. A triangle is regular if its vertices are all regular, and irregular otherwise. Using these definition, we can define the mapping from the barycentric coordinates on a given patch, $\epsilon$, to the limit surface by a weighted average of the effective basis functions associated to \reviewerone{its} 1-ring \cite{stam1998evaluation}. As a result, we can define the limit surface as

\begin{equation}
    \label{eq:Sexp}
        \Gamma(\textbf{r}) = \sum_{i=1}^{N_v} \textbf{c}_{i}\xi_{i}(\textbf{r}),
\end{equation}
where $\mathbf{c}_i$ are vertex locations of the $N_v$ control points, and $\xi_{i}$ is the effective basis function that is associated with quantities associated with $\mathbf{c}_i$ and \reviewerone{has a support $\Gamma_i$}\reviewerone{; note, $\cup_i \Gamma_i = \Gamma$.} The basis functions $\xi_{i}$ span a IGA finite dimensional space $\Psi$ that is the subspace of the Sobolev space $H^{2}(\Gamma)$ \cite{IGA_JCP,Rumpf_IGA}.

To define isogeometric basis sets, we assume that there exists a net of control function values, coincident with the location of the control net. Thus, any scalar function ($f (\mathbf{r})$) can then be expressed in terms of the Loop subdivi  sion basis set via 
\begin{equation}
    f(\mathbf{r})= \sum_{i=1}^{N_v} a_{i}\xi_{i}(\mathbf{r}),
    \label{eq:scalar_func}
\end{equation}
where  $N_{v}$ and $\xi_{i}(\mathbf{r})$ retain the same definition as those prescribed above. The properties of this representation follow from those for subdivision.

Henceforth, the functions $\xi_{i}(\mathbf{r})$ will be referred to as Loop basis. Its properties are (a) positivity, (b) compact support, (c) forming a partition of unity and (d) $C^{2}$ continuity almost everywhere. These properties are critical to the development of both isogeometric analysis as well as defining finite element spaces on the manifold to obtain MHBs.

\section{Current Representation} 

The Loop basis  used to define the geometry, provide the means to define the current as well. To do so, we begin by representing currents on any closed surface $\Gamma$, via the the Helmholtz decomposition as
\begin{equation}
    \textbf{J}(\textbf{r}) = \nabla_{\Gamma}\phi(\textbf{r}) + \nabla_{\Gamma} \cross (\hat{\textbf{n}}\psi(\textbf{r})) + \Bar{\omega}(\textbf{r}),
    \label{eq:HHD}
\end{equation}
where $\Bar{\omega}(\textbf{r})$ is the harmonic field, $\nabla_{\Gamma}$ is  the  surface  gradient, and $\psi(\textbf{r})$ and $\phi(\textbf{r})$ are scalar potentials (that satisfy the mean zero constraint). Assuming that $\Gamma$ is simply connected, $\Bar{\omega}(\textbf{r}) = 0$. While it is possible to develop div-conforming subdivision basis \cite{multiply_subd}, we have chosen to restrict ourselves to simply connected objects\reviewerone{; the rationale being that our basis set offers an exact Helmholtz decomposition, enabling us to develop manifold harmonics for simply connected surfaces that are $C^2$}. In what follows, we construct currents in terms of the scalar potentials using both the loop subdivision basis sets and manifold harmonics.

\subsection{Loop subdivison basis sets}

Using (\ref{eq:scalar_func}) we can define the scalar potentials $\phi (\vb{r})$ and $\psi (\vb{r})$ on the limit surface as
\reviewerone{
\begin{equation}
\begin{split}
     \phi(\textbf{r}) &\approx \tilde{\phi}(\textbf{r}) = \sum_{n=1}^{N_{v}}a^{1}_{n}\xi_{n}(\textbf{r}), \\  
    \psi(\textbf{r}) &\approx \tilde{\psi}(\textbf{r}) = \sum_{n=1}^{N_{v}}a^{2}_{n}\xi_{n}(\textbf{r}).
\end{split} 
\label{eq:f_approx}
\end{equation}
}
It follows from (\ref{eq:HHD}) that it is possible to define the approximation of the current on a simply-connected limit surface as
\begin{subequations}
\label{eq:J_Loop}
\begin{align}
    \begin{split}
    &\textbf{J}(\textbf{r}) \approx \textbf{J}_{N}(\textbf{r}) = \sum_{n}\left[ a^{1}_{n}\textbf{J}_{n}^{1}(\textbf{r}) + a^{2}_{n}\textbf{J}_{n}^{2}(\textbf{r})\right], 
    \end{split}\\
    \begin{split}
    \label{eq:Loop_HH}
    &\textbf{J}_{n}^{1}(\textbf{r}) = \nabla_{\Gamma}\xi_{n}(\textbf{r}),\\
    &\textbf{J}_{n}^{2}(\textbf{r}) = \hat{\textbf{n}}(\textbf{r}) \times \nabla_{\Gamma}\xi_{n}(\textbf{r}).
    \end{split}
\end{align}
\end{subequations}
Finally, since the representation is constructed using conditions on currents that rely on derivatives of the potentials $\tilde{\phi}(\textbf{r})$ and $\tilde{\psi}(\textbf{r})$, leading to the existence of nontrivial solutions to (\ref{eq:J_Loop}), we must enforce uniqueness. In order to ensure uniqueness, we impose an additional zero-mean constraint on the finite dimensional space $\Psi$, leading to 
\begin{equation}
    \Psi = H^{2}(\Gamma) \cap \left \{\int_{\Gamma}f(\textbf{r})d\textbf{r}=0 \right \}.
\end{equation}
A more thorough explanation, as well as, several properties of the basis functions can be found in \cite{Jie,Fu2017GeneralizedDS}.

\subsection{Manifold Harmonics}
While the loop subdivision basis sets are local basis sets, what we explore next is the possibility of a developing a global representation for the potentials $\phi (\vb{r})$ and $\psi(\vb{r})$. In effect, we are seeking the smoothest possible way to interpolate $\psi(\vb{r})$ and $\phi(\vb{r})$; it is well known that the Laplace-Beltrami operator (LBO) is an ideal candidate \cite{Italian_LBO_book,Rosenberg_LBO}. Consider a real-valued function $\chi (\vb{r})$ defined on a compact 2D Riemannian manifold $\Gamma$ embedded in $\mathbb{R}^{3}$. The Laplace-Beltrami operator $\Delta_{\Gamma}$ is defined by
\begin{equation}
    \Delta_{\Gamma} \chi(\vb{r}) := \nabla \cdot (\nabla \chi(\vb{r}) ).
\end{equation}
The LBO $\Delta_{\Gamma}$ admits a complete and countable sequence of eigenfunctions which form an orthonormal basis in $L_{2}\left(\Gamma\right)$ \cite{Italian_LBO_book}, denoted by $\{H_{m}\}$ such that
\begin{equation}
\label{eq:LBOequation}
-\Delta_{\Gamma}H_{m} = \lambda_{m}H_{m}.
\end{equation}
These eigenfunctions, known as Manifold Harmonic Basis (MHB), are the building block for a complete system of eigenfunctions of the vector Laplace–Beltrami operator (or Hodge Laplace operator) $\Vec{\Delta}_{\Gamma} = \nabla_{\Gamma}\mbox{div}_{\Gamma} - \mbox{curl}_{\Gamma} \mbox{curl}_{\Gamma}$. Indeed, the system $\{\nabla_{\Gamma}H_{m},\mbox{curl}_{\Gamma}H_{m}\}$ forms a system of orthogonal nontrival eigenvectors for $\vec{\Delta}_{\Gamma}$ with the same eigenvalues $\lambda_{m}$

\begin{equation}
    -\vec{\Delta}_{\Gamma} \nabla_{\Gamma} H_{m} = \lambda_{m}\nabla_{\Gamma}H_{m},
\end{equation}
\begin{equation}
    -\vec{\Delta}_{\Gamma} \mbox{curl}_{\Gamma} H_{m} = \lambda_{m}\mbox{curl}_{\Gamma}H_{m}.
\end{equation}

Therefore, given $\textbf{J} \in L^{2}(\Gamma)$, we have
\begin{equation}
    \textbf{J} = \sum_{m=1}^{\infty} v_{m}\frac{\nabla_{\Gamma}H_{m}}{\sqrt{\lambda_{m}}} + w_{m}\frac{\mbox{curl}_{\Gamma}H_{m}}{\sqrt{\lambda_{m}}},
\end{equation}
so that $\{\nabla_{\Gamma}H_{m},\mbox{curl}_{\Gamma}H_{m}\}$ is an orthonormal basis for the space of square integrable tangential vector field. 

\begin{figure*}[!t]
\hspace*{-1.4cm}                                                           
\centering
\subfloat[$H_{1}$]{
    \label{fig:boeing_1}
	\includegraphics[width=0.27\textwidth]{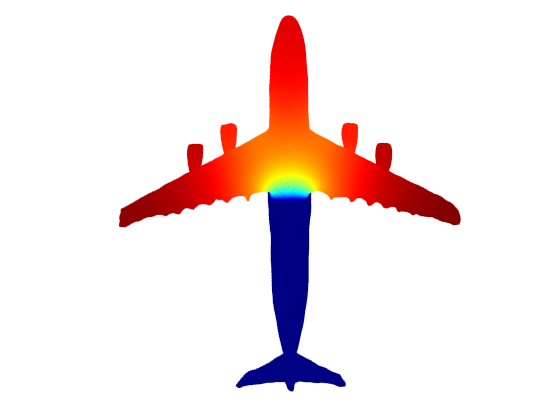} }
	\subfloat[$H_{9}$]{
    \label{fig:boeing_9}
	\includegraphics[width=0.27\textwidth]{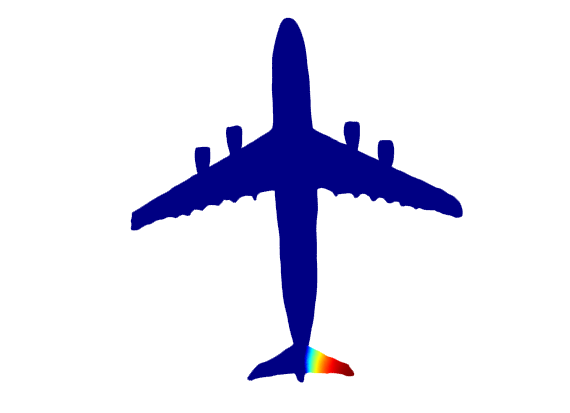} }
\subfloat[$H_{500}$]{
    \label{fig:boeing_500}
	\includegraphics[width=0.27\textwidth]{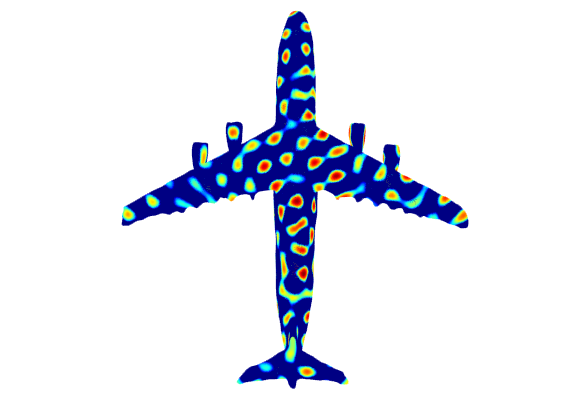} } 
%
\subfloat[$H_{1500}$]{
    \label{fig:boeing_1500}
	\includegraphics[width=0.27\textwidth]{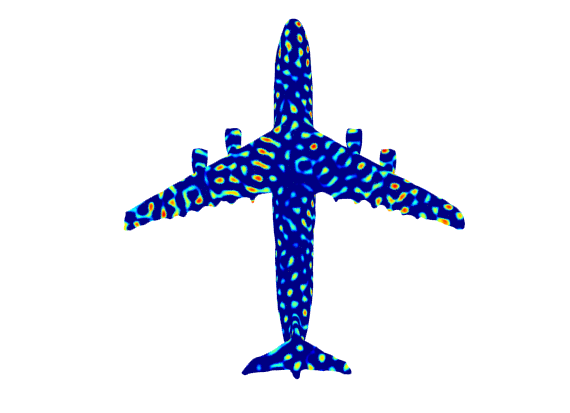} } 
 \\
\caption{ A select few MHs of the jet airliner. (a) $H_{1}$. (b) $H_{9}$. (c) $H_{500}$. (d) $H_{1500}$. }
\label{fig:boeing_MHs}
\end{figure*}

\subsection{Computing the Manifold Harmonics}

In order to numerically compute MHBs, we employ the Loop Subdivision FEM Galerkin method. This is akin to similar efforts using Lagrangian surface descriptions \cite{DLBO,LBO_DNA} that have shown both $h-$ and $p-$ convergence \cite{DLBO,Fu2017GeneralizedDS,LBO_DNA}. The numerics necessary for computing eigenfuctions of the LBO relies on casting the Laplacian eigenvalue problem in a variational setting. The solution of this variational problem is approximated using the finite element Galerkin technique on the surface. We begin by evaluating an inner product of (\ref{eq:LBOequation}) with some test function $v (\mathbf{r}) \in \{\xi_{i}(\mathbf{r})\}$  and then use Green's theorems to arrive to the following:
\begin{equation}
\label{eq:weak_LBO}
\left \langle \nabla_s v (\mathbf{r}), \nabla_s H_{m}(\mathbf{r}) \right \rangle_{\Gamma}  = -\lambda_{m} \left \langle v(\mathbf{r}), H_{m} (\mathbf{r}) \right \rangle_{\Gamma}.
\end{equation}
where $\left \langle f(\textbf{r}), g(\textbf{r}) \right \rangle_{\Gamma} = \int_{\Gamma} f(\textbf{r}) \cdot g(\textbf{r}) d\textbf{r}$ follows the standard inner product definition. The MH $H_m (\vb{r})$ is represented in the same fashion as (\ref{eq:scalar_func}) leading to 
\begin{equation}
    \label{eq:MH_loop}
    H_{m} \approx \widetilde{H}_{m}(\mathbf{r}) = \sum_i^{N_v} h^{i}_m \xi_i (\mathbf{r}),
\end{equation}
for $h^{i}_m \in \mathbb{R}$. This leads to a  generalized eigenvalue problem 
\begin{equation}
    [A] [H] = -[\Lambda] [B][H],
    \label{eq:geneigvalprob}
\end{equation}
where,
\begin{subequations}
\begin{equation}
    \begin{aligned}
     \left[A\right]_{ij} = \int_{\Gamma_{i}} \nabla_s \xi_{i} (\mathbf{r}) \cdot \nabla_s \xi_{j} (\mathbf{r}) d \mathbf{r},
    \end{aligned}
    \label{eq:stiffmatrix}
\end{equation}
\begin{equation}
    \begin{aligned}
    \left[B\right]_{ij} = \int_{\Gamma_{i}} \xi_{i} (\mathbf{r})  \xi_{j} (\mathbf{r}) d\mathbf{r}.
    \end{aligned}
    \label{eq:massmatrix}
\end{equation}
\end{subequations}
For this generalized symmetric eigenvalue problem  $[A] \in \mathbb{R}^{N_v \times N_v}$ is positive semi-definite, $[B] \in \mathbb{R}^{N_v \times N_v}$ is positive definite, $[\Lambda] \in \mathbb{R}^{N_v \times N_v}$ contains $N_v$ eigenvalues along its diagonal, and $[H] \in \mathbb{R}^{N_v \times N_v}$ contains the solution vectors, i.e. the coefficients of each eigenvector defined in (\ref{eq:MH_loop}), in \reviewerone{its} column space. For this symmetric generalized eigenvalue problem we have $[H]^{T}[A][H] = [\Lambda]$ and $[H]^{T}[B][H] = [I]$, where $[I]$ is the identity matrix. From the previous relations, it follows that the eigenfunctions are orthogonal with respect to the $[B]$-based scalar product (i.e., $\langle \textbf{H}_{i}, \textbf{H}_{j} \rangle_{[B]}  = \textbf{H}_{i}^{T}[B]\textbf{H}_{j}$). \reviewerone{The eigenvectors with corresponding eigenvalues can be calculated with a direct eigensolver or by using the efficient band-by-band computation method presented in \cite{MHs}. There is a extensive body of literature on efficient computation of these functions, largely applied to computer graphics \cite{LBO_Zhang}.}

Given the representation of each of the eigenfunction, it follows that functions defined on the manifold can be written in terms of these eigenfunctions, as can its divergence and \reviewerone{rotational}. Specifically, 

\reviewerone{
\begin{subequations}
\begin{equation}
    \begin{aligned}
     \nabla_{\Gamma} H_m(\mathbf{r}) \approx \nabla_{\Gamma} \widetilde{H}_{m}(\mathbf{r})  = \sum_i^{N_v} h^{i}_m\nabla_{\Gamma} \xi_i(\mathbf{r}),
    \end{aligned}
\end{equation}
\begin{equation}
    \begin{aligned}
   \mbox{curl}_{\Gamma} H_m (\mathbf{r}) \approx \mbox{curl}_{\Gamma} \widetilde{H}_{m} (\mathbf{r}) = \sum_i^{N_v} h^{i}_m \hat{\textbf{n}} \times \nabla_{\Gamma} \xi_i(\mathbf{r}).
    \end{aligned}
\end{equation}
\end{subequations}
}
Using these expressions, the currents may alternatively be written in in terms of this basis as 

\begin{subequations}
\label{eq:J_MH}
\begin{align}
    \begin{split}
    \label{eq:MH_J}
    &\textbf{J}(\textbf{r}) \approx \textbf{J}_{M}(\textbf{r})  =  \sum_{m=1}^{N_v}\left[ v_{m} \tilde{\textbf{J}}_{m}^{1}(\textbf{r}) + w_{m}\tilde{\textbf{J}}_{m}^{2}(\textbf{r})\right], 
    \end{split}\\
    \begin{split}
    \label{eq:MH_HH}
       &\tilde{\textbf{J}}_{m}^{1}(\textbf{r}) = \frac{\nabla_{\Gamma}{\widetilde{H}_{m}}(\textbf{r})}{\sqrt{\lambda_{m}}},\\
           &\tilde{\textbf{J}}_{m}^{2}(\textbf{r}) = \frac{\mbox{curl}_{\Gamma}{\widetilde{H}_{m}}(\textbf{r})}{\sqrt{\lambda_{m}}}.
    \end{split}
\end{align}
\end{subequations}
such that $\{\textbf{J}_{m}^{1},\textbf{J}_{m}^{2}\}$ is an orthonormal basis of the space of $L^{2}$ tangential vector fields. Similarly as stated above, we impose uniqueness of (\ref{eq:J_MH}) using a zero-mean constraint.

\subsection{Illustration of Manifold Harmonic Transforms}
While the manifold harmonic transform has been commonplace in the compute graphics literature for an array of applications, to the authors' knowledge it has not been utilized in analysis of physics on manifolds.  In particular, one of \reviewerone{its} many attractive features is its ability to rigorously compress the system. In what follows, we illustrate some of the features of this approach within the framework of this paper. To wit, we consider representation of currents on two different objects: a bumpy cube and a jet airliner. Our goal is \textcolor{green!50!black}{to} examine the convergence of the representation of the current to a bandwidth of $M$ harmonics.  

In both instances, we reconstruct a surface current generated by a  1 GHz plane wave incident in the $-\hat{\textbf{z}}$, respectively. In  Fig.~\ref{fig:bumpy_cube_MHs}, we visualize the \reviewerone{manifold harmonic representation} of the current for a bumpy cube, \reviewerone{with 5124 DoFs}, and in Fig.~\ref{fig:boeing_MHs} for a jet airliner, \reviewerone{with 12132 DoFs} . As can be seen in both figures, the first $\textbf{J}_{m}$ functions capture the coarse features of the current and the next, high frequency ones, correspond to the details.
\begin{table}[h!]
\begin{center}

\resizebox{\columnwidth}{!}{
    \begin{tabular}{c c c c c c c}
\hline
\multirow{2}{*}{Bumpy Cube} &  & \multirow{1}{*}{$M$} &200  &1000 & 2000 & 5122 \\
& & $\epsilon$  &9.87E-4  & 4.02E-4 & 4.16E-5 & 3.65E-17 \\

\multirow{2}{*}{Jet airliner} & & \multirow{1}{*}{$M$} &500  &1000 & 2000& 12130 \\
& & $\epsilon$  &4.21E-4  & 9.93E-5 & 4.14E-5 & 9.03E-17 \\
\hline
\end{tabular}
}
\end{center}
\caption{Relative $\epsilon$ error, \reviewerone{with respect to number of manifold harmonics M}, in the reconstructed surface currents density.}
\label{table1}
\end{table}
Table.~\ref{table1} demonstrates the precision of the inverse MHT (\ref{eq:J_MH}) w.r.t original current $\textbf{J}_{N}$ as we increase the number of MHs. Our metric for validation is the reconstruction error  $\epsilon = \norm{\textbf{J}_{N}(\textbf{x}_{i}) - \tilde{\textbf{J}}(\textbf{x}_{i})}_{[B]} ^{1/2} $. Note, $\textbf{J}_{N}(\textbf{x})$ is the current on the surface as approximated by the Loop-subdivision basis set. In both candidate objects, we find that as expected, $\epsilon$ decreases as the number of MH $M$ increases, eventually approaching machine precision.

\begin{figure*}[!t]
\hspace*{-1cm}                                                           
\centering
\subfloat[$|\tilde{\textbf{J}}_{200}|$]{
    \label{fig:bumpy_cube_200}
	\includegraphics[width=0.2\textwidth]{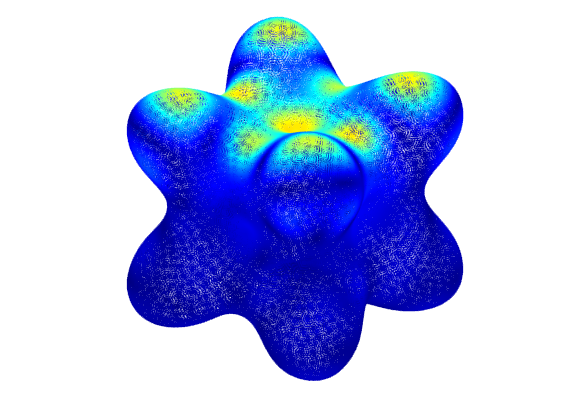} }
	\subfloat[$|\tilde{\textbf{J}}_{1000}|$]{
    \label{fig:bumpy_cube_1000}
	\includegraphics[width=0.2\textwidth]{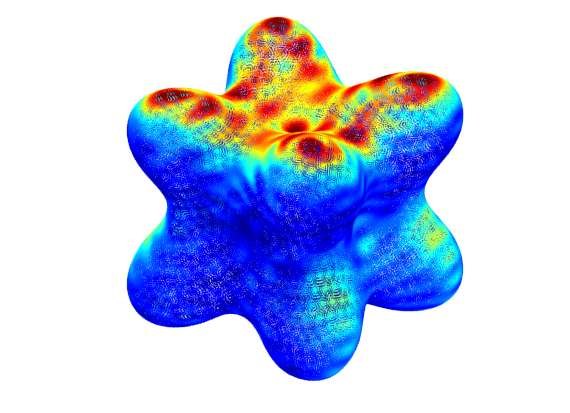} }
\subfloat[$|\tilde{\textbf{J}}_{2000}|$]{
    \label{fig:bumpy_cube_2000}
	\includegraphics[width=0.2\textwidth]{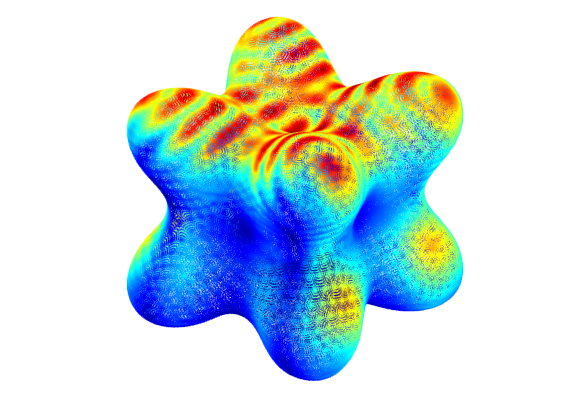} } 
%
\subfloat[$|\tilde{\textbf{J}}_{5122}|$]{
    \label{fig:bumpy_cube_5122}
	\includegraphics[width=0.2\textwidth]{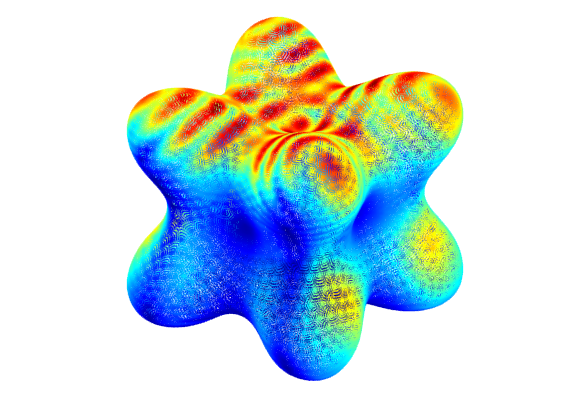} } 
%
\subfloat[$|\textbf{J}|$]{
    \label{fig:bumpy_cube}
	\includegraphics[width=0.2\textwidth]{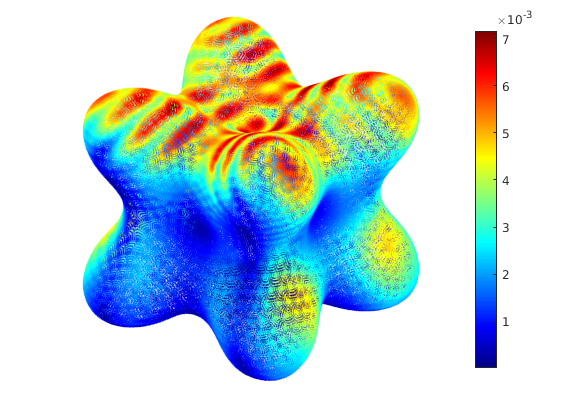} }
 \\
\caption{  Reconstruction of the target current $\textbf{J}$ obtained with an increasing number of MHs for the bumpy cube. }
\end{figure*}

\section{Field Solvers}

Thus far, we have discussed Loop subdivision basis and its mapping to MHBs. In this Section, we detail the discretization of \eqref{eq:CCCFIER}, in terms of \reviewerone{both loop basis and MHB;} in particular, we use a Galerkin prescription to discretize these equations. \reviewerone{In what follows, we describe both. }
\reviewerone{\subsection{Discretization using Loop subdivision basis}}
Note that discretizing Calder\'{o}n type operators requires intermediate spaces, effected through a Gram matrix. We define the required Gram-matrix $[G]$ using
\begin{equation}
    [G]_{nm}^{lk} = \delta_{lk} \left \langle  \textbf{J}_n^l, \textbf{J}_m^k \right \rangle_{\Gamma_{n}}, 
\end{equation}
where $\delta_{lk}$ is a Kronecker's delta  \cite{Jie}. In effect, the system of matrices to be solved can be written as 
\begin{subequations}
\begin{equation}
    \label{eq:matrix_eq}
    \left [ Z\right ] \left [I \right ] = \left [ V\right ]
\end{equation}
where,
 
\begin{equation}
    \left [ Z \right] =[G]^{-1}\left [  \left [ L \right ] + \left [ K \right ] \right ]
\end{equation}
with
\begin{equation}
[K]_{nm}^{lk} = \left \langle \textbf{J}_n^l (\vb{r}), \frac{\textbf{J}^k_m}{2} (\textbf{r}) - \mathcal{K}_{\kappa}\circ \textbf{J}_m^k(\textbf{r}) \right \rangle_{\Gamma_{n}},
\end{equation}
\begin{equation}
    \label{eq:T}
    [T]_{\widetilde{\kappa},nm}^{lk} = \left \langle \textbf{J}_n^l (\vb{r}),  \mathcal{T}_{\widetilde{\kappa}}\circ \textbf{J}_m^k(\textbf{r}) \right \rangle_{\Gamma_{n}},
\end{equation}
where $\widetilde{\kappa} \in \{\kappa',\kappa\}$, and, as defined earlier $ \kappa' = \kappa + 0.4\varsigma^{2/3}\kappa^{1/3}$, and $\varsigma$ is the mean curvature of the object, and, 
\begin{equation}
[L] = -2[T]_{\kappa'}[G]^{-1}[T]_\kappa.
\end{equation}
\end{subequations}
Furthmore, we have

\begin{subequations}
\reviewerone{
\begin{equation}
    [I]_{m}^{k} = a^k_{m},
\end{equation}
}
\begin{equation}
    [V]_{n}^{k} = [G]^{-1}\left[-2[T]_{\kappa'}[G]^{-1}[V_{T}]_{n}^{k} +[V_{K}]_{n}^{k}\right],  
\end{equation}
\end{subequations}
with  
\begin{subequations}
\begin{equation}
    [V_{T}]_{n}^{k}  = \left \langle \textbf{J}_{n}^{k}(\textbf{r}),   \textbf{E}^{i}(\textbf{r})\right \rangle_{\Gamma_{n}},
\end{equation}
\begin{equation}
    [V_{K}]_{n}^{k} = \left \langle \textbf{J}_{n}^{k}(\textbf{r}),\hat{\textbf{n}}  \cross \textbf{H}^{i}(\textbf{r})\right \rangle_{\Gamma_{n}},  
\end{equation}
\end{subequations}
Lastly, we note that the stabilizing properties of the Calderón preconditioner are local \cite{Local_CC_CFIER,MLFMA_CC_CFIER}, which allows the use of a localized version of the preconditioner $[T]_{\kappa'}$. As such, we choose to omit all interactions of a distance greater than \reviewerone{$1.25\lambda$}.  

\subsubsection{Wideband MLFMA for Evaluation of Inner Products}

At this point, we note the following: the domain of support of each \reviewerone{loop} basis function is electrically large and are on average $\approx 0.9\lambda$; this makes efficient evaluation of both inner products \emph{and} matrix vector products challenging. Furthermore, we note that the basis functions are \reviewerone{third} order \reviewerone{and the geometry is fourth order}. Both serve to exacerbate costs \reviewerone{as one needs higher order quadrature rules over both test and source domains}. To ameliorate these, we exploit the wide-band FMM introduced by the authors in \cite{vikram2009novel,dault2016mixed}. The framework we propose has been used to accelerate matrix evaluations as well as matrix vector products for the Generalized Method of Moments (GMM) wherein patch sizes can be several wavelengths long \cite{dault2016mixed} using a mixed potential formulation. 

\reviewerone{In order to motivate this framework we begin by briefly examining the fundamental rubrics involved in evaluating all inner products. As has been shown in \cite{multiply_subd}, expected convergence rates for for the EFIE are obtained, \emph{albeit} with RWG type basis functions, \emph{provided} one uses an adaptive quadrature based on subdividing each parent triangle. This implies that it is possible to choose a refinement level such that one can use a piecewise flat approximation for evaluation of singular integrals. But the downside is clear; the cost of evaluating the inner-products is high. Indeed, if there are $N_p$ quadrature points per basis, the cost for evaluation of of matrix interactions scale as $\mathcal{O}(4 N_p^2 N_v^2)$. Our goal is to reduce this to $\mathcal{O}(2N_{p}N_{v}\log\left(N_{p}N_{v}\right))$.}

\reviewerone{ In order to mitigate the high computational cost we employ MLFMA. We note that the following is true for all methods that use higher order modeling. To set the stage for the discussion, consider a matrix element $[T]_{\tilde{\kappa},nm}^{lk}$. It comprises contributions from both the magnetic vector potential and electric scalar potential. Let us focus on the latter, specifically, just evaluating the scalar potential due to $\vb{J}_n^1$,   
\begin{equation}\label{eq:potMethod1}
    \Phi (\vb{r}) = \frac{-j}{\omega\varepsilon}\int_{\Gamma_n}\nabla' G_{\tilde{\kappa}} (\vb{r}, \vb{r}')  \cdot \vb{J}^{1}_{n} (\vb{r}')d\vb{r}' .
\end{equation}
}
\reviewerone{Consider a 1-level MLFMA prescription for an alternative evaluation of \eqref{eq:potMethod1}. We denote the center of a leaf box by $\vb{r}_s$; at observations points sufficiently far away, where $\vb{X} = \vb{r} - \vb{r}_{s}$, the potential $\Phi (\vb{r})$ can also be evaluated using 
\begin{subequations}\label{eq:potFMM}
\begin{equation}
\Phi (\vb{r}) = \frac{- j\tilde{\kappa}^2 }{16 \pi^2 \omega \varepsilon}\int_{S^2} \hat{\vb{k}} \cdot \mathcal{M} \left ( \tilde{\kappa} \hat{\vb{k}}, \vb{r}_s \right ) \mathcal{T} \left (\tilde{\kappa}\hat{\vb{k}}, \vb{X} \right ) d^2 \hat{\vb{k}},
\end{equation}
where 
\begin{equation}\label{eq:potMuyltipole}
\mathcal{M} \left ( \tilde{\kappa}\hat{\vb{k}}, \vb{r}_s \right ) = \int_{\Gamma_n} \vb{J}^{1}_{n} (\vb{r}') e ^{- j \tilde{\kappa}\hat{\vb{k}} \cdot (\vb{r}_s - \vb{r}') } d\vb{r}', 
\end{equation}
is the source to multipole map, and 
\begin{equation}
    \mathcal{T}\left (\tilde{\kappa}\hat{\vb{k}}, \vb{X} \right ) \doteq \sum_{n=0}^{\infty} (-j)^n (2n+1) h_n^{(2)}(\tilde{\kappa}X) P_n\left(\hat{\vb{k}}\cdot\hat{\vb{X}}\right),
\end{equation}
\end{subequations}
is the translation operator. Here, $S^2$ denotes the unit sphere, parametrized by $(\theta,\phi)\in\left[0,\pi\right]\times\left[0,2\pi\right]$. We note that $\hat{\vb{k}}=\hat{\vb{k}}(\theta,\phi)$. Note, the gradient on the Green's function is evaluated spectrally; furthermore, a traditional approach to using MLFMA would ensure that the entire support domain $\Gamma_n$ lies within a leaf box, i.e., for loop basis functions the size of the leaf box is $\Delta_0 \approx 0.9\lambda$. This means that that each leaf box has approximately $N_{p}$ quadrature points, and the cost of computing the near field interactions in the MLFMA scheme is $\mathcal{O}(4N_{p}^2N^{2}_{v})$. As $N_p$ is relatively high, this still untenable. We would like $\Delta_0$ to be as small as possible, such that it contains far fewer quadrature points.}

\begin{figure}[!htbp]
  \centering
  \includegraphics[width=0.25\textwidth]{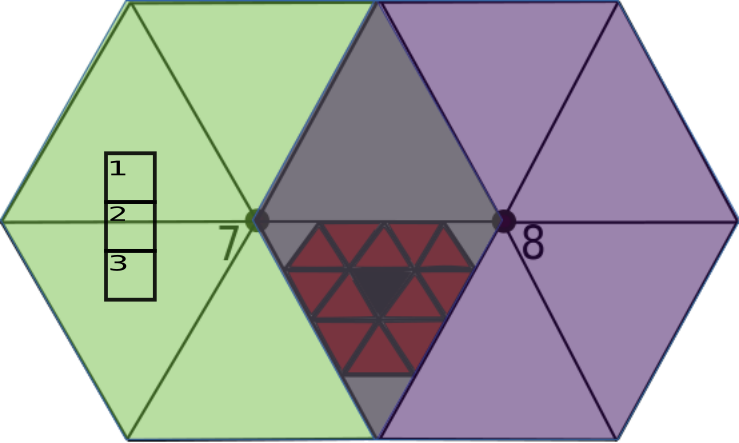}
  \caption{\reviewerone{The support of two basis (no. 7 and no. 8) is shown in green and purple, respectively. In computing the interaction between the two, patch (triangular tesselation of control nodes) is partitioned into sub-patches (shown in red) to create an adaptive quadrature. The entire object is embedded in an MLFMA tree, and the size of a leaf box (shown within the green region) is about the size of a sub-patch.}}
  \label{fig:basis}
\end{figure}

\reviewerone{Consider instead Fig. \ref{fig:basis} which shows three leaf level boxes within $\Gamma_n$. Furthermore, assume that we have to compute the self interaction of basis $n$ using MLFMA, boxes (1,3) are in the far field of each other, and (2,3) and (1,2) are in the near field of each other. Lets re-examine the evaluation of $\Phi(\vb{r})$ for $\vb{r} \in \Gamma_{3,n}$. Using \eqref{eq:potMethod1}
\begin{equation}\label{eq:phi3}
\begin{split}
    \Phi (\vb{r}) & = \frac{-j}{\omega\varepsilon}\int_{\Gamma_{1,n}} \nabla' G_\kappa (\vb{r}, \vb{r}')  \cdot \vb{J}^{1}_{n} (\vb{r}')d\vb{r}' \\
    & - \frac{j}{\omega\varepsilon}\int_{\Gamma_{2,n}} \nabla' G_\kappa (\vb{r}, \vb{r}') \cdot \vb{J}^{1}_{n} (\vb{r}')d\vb{r}', \\
    & = \Phi_{13} (\vb{r}) + \Phi_{23}(\vb{r}).
\end{split}
\end{equation}  
Here, $\Gamma_{i,n}$ denotes the intersection of box-$i$ with $\Gamma_n$. There are two possible ways of evaluating the farfield interaction,  $\Phi_{13}(\vb{r})$, using a variation of \eqref{eq:potFMM} with the understanding that the the domain of integration in \eqref{eq:potMuyltipole} is confined to $\Gamma_{1,n}$. Specifically, 
\begin{subequations}
\begin{equation}
\begin{split}
    \Phi_{13} (\vb{r}) & = \frac{- j\tilde{\kappa}^2 }{16 \pi^2 \omega \varepsilon}\int_{S^2} \hat{\vb{k}} \cdot \mathcal{M}_1 \left ( \tilde{\kappa} \hat{\vb{k}}, \vb{r}_s \right ) \mathcal{T} \left (\tilde{\kappa}\hat{\vb{k}}, \vb{r} - \vb{r}_s \right )d^2 \hat{\vb{k}}\\
    & = \frac{\tilde{\kappa} }{16 \pi^2 \omega \varepsilon}\int_{S^2}  \mathcal{M}_2 \left ( \tilde{\kappa} \hat{\vb{k}}, \vb{r}_s \right ) \mathcal{T} \left (\tilde{\kappa}\hat{\vb{k}}, \vb{r} - \vb{r}_s \right )d^2 \hat{\vb{k}}
\end{split}
\end{equation}
where 
\begin{equation}\label{eq:potMLFMASpectral}
\mathcal{M}_1 \left ( \tilde{\kappa}\hat{\vb{k}}, \vb{r}_s \right ) = \int_{\Gamma_{1,n}} \vb{J}^{1}_{n} (\vb{r}') e ^{- j \tilde{\kappa}\hat{k} \cdot (\vb{r}_s - \vb{r}') } d\vb{r}', 
\end{equation}
\begin{equation}\label{eq:potMLFMASpatial}
\begin{split}
    \mathcal{M}_2 \left ( \tilde{\kappa}\hat{\vb{k}}, \vb{r}_s \right )   & = -\int_{\Gamma_{1,n}} \nabla'\cdot \vb{J}^{1}_{n} (\vb{r}') e ^{- j \tilde{\kappa}\hat{\vb{k}} \cdot (\vb{r}_s - \vb{r}') } d\vb{r}' \\
    & + \int_{\partial \Gamma_{1,n}} \hat{\vb{u}}_{\Gamma_{1,n}} (\vb{r}') \cdot \vb{J}^{1}_{n} (\vb{r}') e ^{- j \tilde{\kappa}\hat{\vb{k}} \cdot (\vb{r}_s - \vb{r}') } d\vb{r}',
\end{split}
\end{equation}
\end{subequations}
and $\hat{\vb{u}}_{\partial \Gamma}$ is outward pointing normal to the boundary $\partial \Gamma$. Consider next, the near field evaluation of $\Phi_{23} (\vb{r})$. As we want the leaf box size to be small, the minimum distance between the box centers becomes very small, and as such the order of singularity, due to the gradient on the Green's function, introduces near-singular integration challenges. The remedy that is typically taken is to transfer the derivative onto the basis function. Specifically, 
\begin{equation}\label{eq:potMethod2Alt}
    \begin{split}
            \Phi_{23}(\vb{r}) & =  \frac{j}{\omega\varepsilon}\int_{\partial \Gamma_{2,n}}  G_\kappa (\vb{r}, \vb{r}') \hat{\vb{u}}_{\partial \Gamma_{2,n}}(\vb{r}') \cdot \vb{J}^1_n (\vb{r}')d\vb{r}' \\
    & - \frac{j}{\omega\varepsilon}\int_{\Gamma_{2,n}}  G_\kappa (\vb{r}, \vb{r}') \nabla'\cdot \vb{J}^1_n (\vb{r}')d\vb{r}'.
    \end{split}
\end{equation}
The lessons we take from the above equations are as follows: (a) the aforementioned line integrals have to be accounted for in \eqref{eq:potMethod2Alt} as they are implicitly included  in \eqref{eq:potMLFMASpectral} and should cancel on the shared interface; (b) unfortunately, finding intersections between higher order surfaces and boxes is non-trivial; (c) with all challenges considered, we have to use \eqref{eq:potMethod2Alt} to evaluate $\Phi_{23}(\vb{r})$ and \emph{not} \eqref{eq:phi3}; Alternatively, it can be proven that interior line integral should vanish. This implies that an ideal choice would be to use \eqref{eq:potMLFMASpatial} and \eqref{eq:potMethod2Alt} \emph{sans} line integrals in both to evaluate $\Phi_{13}(\vb{r})$ and  $\Phi_{23}(\vb{r})$, respectively. Note, this example is illustrative. Further complication arises in the evaluation of the electric field, as it calls for the gradient of the scalar potential. As a result, one needs additional line integrals to reduce the singularity. As is evident from the above discussion, using a mixed potential formulation together with wideband MLFMA permits evaluation of all integrals, near and far, without the consideration of the troublesome line integrals, but at the cost of more tree traversals \cite{vikram2009novel}. Indeed, the size of the leaf box can now be as small as computationally expedient. Leaf box sizes can be chosen such that it contains $\mathcal{O}(1)$ quadrature points, reducing the cost of near field evaluation to $\mathcal{O}(2 N_p N_{v})$. }

\reviewerone{ We elucidate this process by applying it to equation
(\ref{eq:T}), as depicted in Fig.~\ref{fig:basis}; (\ref{eq:T}) contains four independent terms that must be computed in the inner integral: three corresponding to the vector potential, and one corresponding to the scalar potential. It follows that any matrix element can be computed in terms of its partial contributions such that
\begin{equation}\label{eq:innerProduct}
    \begin{split}
        [T]_{\widetilde{\kappa},nm}^{lk} = \sum_\zeta \sum_\gamma -j\omega \mu \left \langle \vb{J}^{l}_{n},  \hat{\vb{n}}\left( \vb{r} \right) \times \mathcal{S}^{\gamma}_{\widetilde{\kappa}} \circ \vb{J}_{m}^{k} \right \rangle_{\Gamma_{n}^{\zeta}} \\+ \frac{j\delta_{l2}\delta_{1k}}{\omega \epsilon}  \left \langle \nabla_{\Gamma} \cdot \vb{J}_{n}^{1},  \mathcal{S}^{\gamma}_{\widetilde{\kappa}} \circ \nabla_{\Gamma} \cdot \vb{J}_{m}^{1} \right \rangle_{\Gamma_{n}^{\zeta}} ,
    \end{split}
\end{equation}
where,
\begin{equation}\label{eq:innerProduct}
     \mathcal{S}^{\gamma}_{\widetilde{\kappa}} \circ \vb{J}_{m}^{k}  = \int_{\Gamma_{m}^{\gamma}} G_{\tilde{\kappa}} (\vb{r}, \vb{r}') \vb{J}^{k}_{m}\left(\vb{r}'\right) d\vb{r}', 
\end{equation}
and the indices $\zeta$ and $\gamma$ are subpatches of $\Gamma_m$ and $\Gamma_n$; subpatches within each other's farfield are constructed via MLFMA whereas nearfield patches are constructed via direct integration, see Fig.~\ref{fig:basis}.}

\begin{figure*}[!t]
\hspace*{-1.8cm}                                                
\centering
\subfloat[$|\tilde{\textbf{J}}_{500}|$]{
    \label{fig:boeing_500_J}
	\includegraphics[width=0.18\textwidth]{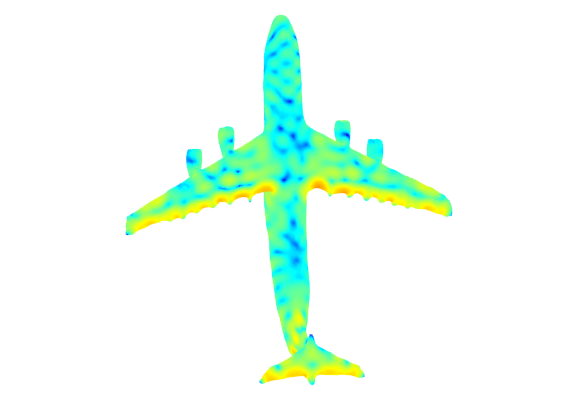} }
	\subfloat[$|\tilde{\textbf{J}}_{1000}|$]{
    \label{fig:boeing_1000}
	\includegraphics[width=0.18\textwidth]{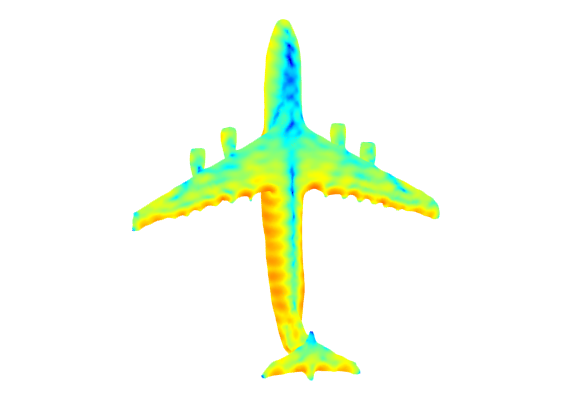} }
\subfloat[$|\tilde{\textbf{J}}_{2000}|$]{
    \label{fig:boeing_2000}
	\includegraphics[width=0.18\textwidth]{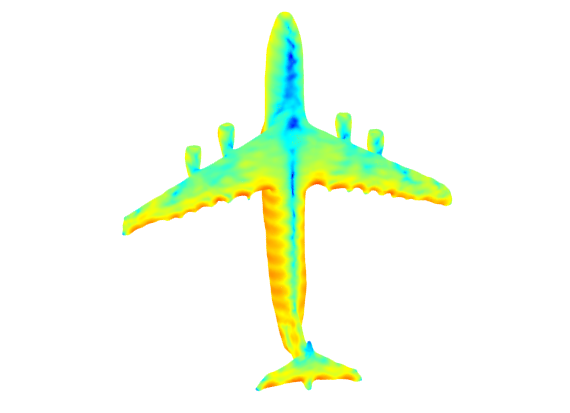} } 
%
\subfloat[$|\tilde{\textbf{J}}_{12130}|$]{
    \label{fig:boeing_12130}
	\includegraphics[width=0.18\textwidth]{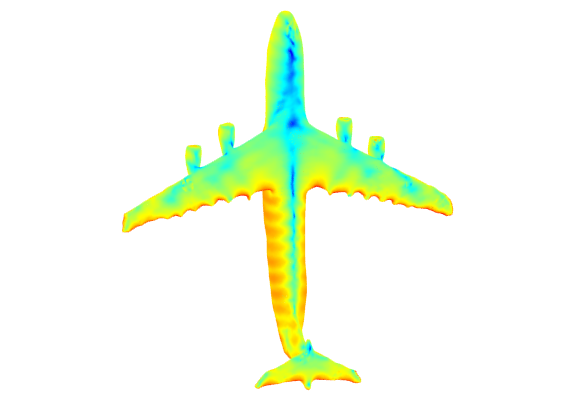} } 
%
\subfloat[$|\textbf{J}|$]{
    \label{fig:boeing}
	\includegraphics[width=0.18\textwidth]{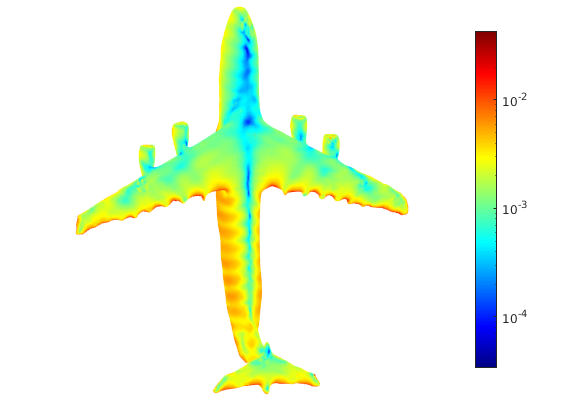} }
 \\
\caption{ Reconstruction of the target current $\textbf{J}$ obtained with an increasing number of MHs for the jet airliner.}
\end{figure*} 

\begin{figure}[!ht]
  \centering
  \includegraphics[width=0.35\textwidth]{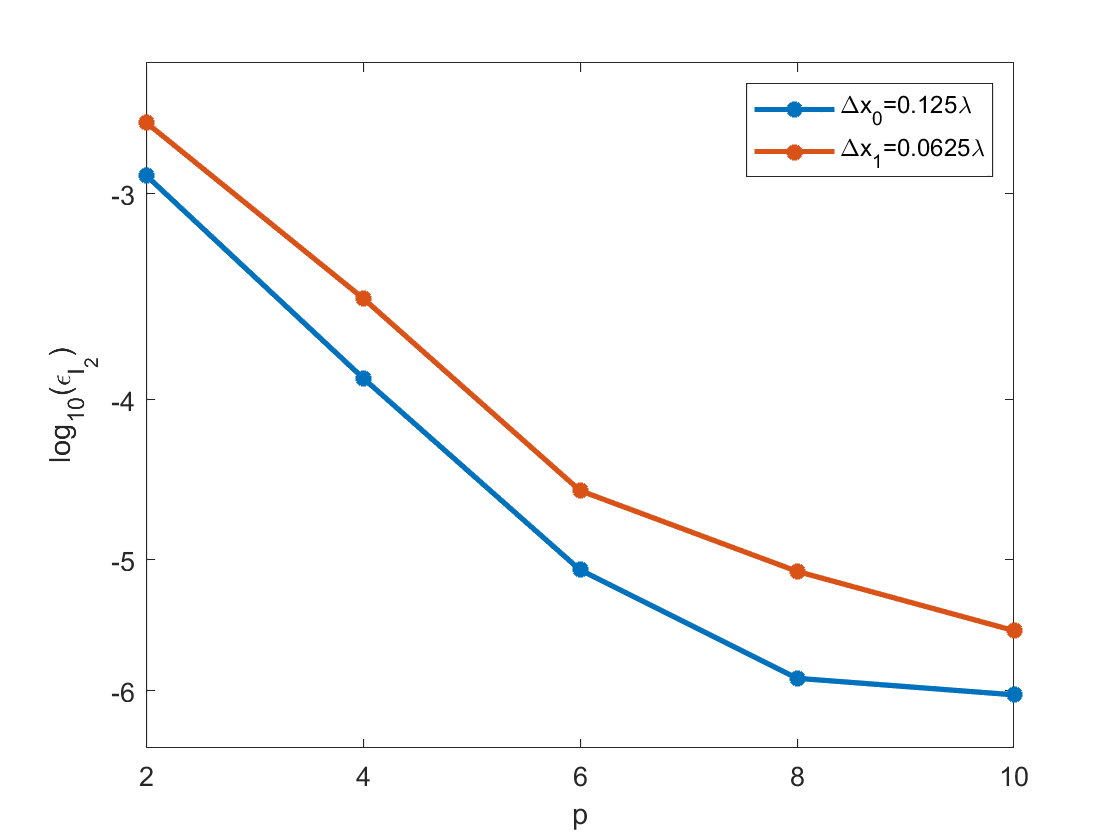}
  \caption{\reviewerone{ Convergence in relative $l_{2}$ error of partial matrix element using wideband MLFMA versus the direct fill algorithm, where $p$ is the expansion coefficient} }
  \label{fig:mlfma_near}
\end{figure} 

\subsection{Manifold Harmonic Transform of CC-CFIER}  

As presented above, the MHs are constructed as a linear combination of the loop subdivision basis functions, and can thus be seen as global basis functions built on top of loop subdivision basis set. As an aside, these basis would be excellent candidate to create a reduced order representation of currents. Consider a reduced $M$ orthogonal MHBs that span $W^{MH}(\Gamma) \subset \Psi(\Gamma)$. This is tantamount to using $M<n_v$ for both the representation and measurement space in \eqref{eq:J_MH}. As a result, one obtains a compressed impedance matrix. \reviewerone{As such we can reconstruct (\ref{eq:matrix_eq}) as}

\reviewerone{
\begin{equation}
    [H][Z_{H}][I_{H}] = [H][V_{H}]
\end{equation}
}
\reviewerone{where,}
\begin{subequations}
\reviewerone{
\begin{equation}
    \left [ Z_{H} \right] = \left [ L_{H} \right ] + \left [ K_{H} \right ],
\end{equation}
\begin{equation}
    \left [ I_{H} \right] = \left [H\right]^{T} \left[ I \right],
\end{equation}
\begin{equation}
    [V_{H}] = \left[H\right]^{T}\left[-2[T]_{\kappa'}\left [H\right]\left[H\right]^{T}[V_{T}] +[V_{K}]\right],  
\end{equation}
\mbox{  and,}
\begin{equation}
[L_{H}] = -2[H]^{T}[T_{\kappa'}][H][H]^{T}[T_\kappa][H],
\end{equation}
\begin{equation}
[K_{H}] = -2[H]^{T}[K][H].
\end{equation}
}
\end{subequations}

\reviewerone{Here $\left [H\right]^{T}$ denotes the transpose of $\left [H\right]$.}
    
\section{Numerical Examples}
In this section, we present a collection of numerical results to demonstrate the efficacy of the proposed approach. \reviewerone{All examples in this section were run serially using a single 2.4 GHz Intel Xeon Gold 6148 CPU core on the HPC Center at Michigan State University. As 
parallelization is \emph{not} used to generate neither the near field matrix elements or MATVECs, RWG data is only generated for the number of DoFs that is within reach of the available resources. Note, the number of degrees of freedom necessary to model the structure using piecewise flat triangles (and indirectly, the number of RWG basis) is significantly higher. } As alluded to in the introduction, the two main contributions are (a) \reviewerone{a fast method for evaluating matrix vector products for higher order geometries and higher order physical basis} illustrated via application to subdivision based isogeometric formulation for simply connected objects, and (b) employing manifold harmonics for EM analysis. To this end, the data presented in this section highlights the following:  (i) the accuracy of the two proposed approaches when compared against analytical data; (ii) the improved spectral properties of the CC-CFIER by means of the reduced numbers of iterations required for convergence of the GMRES iterative solver for Loop and MHB, (iii) the high-accuracy and reduced DOF under the MHB, and (iv) application of both to analyzing complex targets.

Unless otherwise stated, we compute scattering due to a  plane wave field propagating in $\hat{\kappa} = -\hat{z}$ and polarized along $\hat{x}$ axis. Furthermore, we compare radar cross sections (RCS) in the $\phi = 0$ plane, using the proposed methods against either analytical data or a validated method of moments code that is based on RWG basis functions, otherwise referred to as RWG-CFIE. For every scattering experiment presented in the tables, the maximum relative far-field error, denoted by $\epsilon_{\infty}$, is defined as 
\begin{equation}
    \epsilon^{\infty} = \frac{\max\limits_{\hat{\textbf{x}}}|\textbf{E}_{\infty}^{calc}(\hat{\textbf{x}}) - \textbf{E}_{\infty}^{ref}(\hat{\textbf{x}})|
    }{\max\limits_{\hat{\textbf{x}}}|\textbf{E}_{\infty}^{ref}(\hat{\textbf{x}})|},
\end{equation}
where the reference solutions $\textbf{E}_{\infty}^{ref}$ was computed by Mie series in the case of spherical scatterers, otherwise, by a Loop subdivision based CC-CFIER. All of the numerical results presented in the tables and graphs in this section were obtained by prescribing a GMRES residual tolerance equal to $10^{-5}$ for the overall system and $10^{-11}$ for inverting the gram matrix with a diagonal preconditioner. Finally, we note that we provide the iteration count to reach the specified GMRES tolerance, the time taken to reach the prescribed tolerance, and the error relative to the benchmark data. 

\reviewerone{\subsection{Accuracy of wideband MLFMA for adaptive interactions} Herein, we study the the accuracy of using wideband MLFMA is to alleviate the computational complexity associated with nearfield computations. To test the controllable accuracy of the aforementioned scheme, we conduct a controlled test. As an aside, the support of a basis function is a one-ring associated with a control vertex, three basis functions are defined on a patch. The most efficient assembly of interactions is computing these in a patchwise manner. To that end, consider two patches that share an edge. The edge length of each patch is approximately $0.25\lambda$. We compute the patch to patch interaction by using an integration rule developed by subdividing each patch into 16 sub-patches and using a 3-point rule in each. Note, we are not computing self-patch interactions. Next, we compute the same interaction, but through an MLFMA tree with leaf box sizes $\Delta_0 = 0.125 \lambda$ and $\Delta_0 = 0.0625 \lambda$ which results in 1-level and 2-level computation of the interactions. The standard tree partitioning of interactions is used; the leaf box of size $0.125\lambda$ has about 4 subpatches, whereas $0.0625\lambda$ has approximately 1 subpatch. Given the size of leaf box, interactions are computed using wideband MLFMA which invokes Accelerated Cartesian Expansions (ACE) for leaf box sizes smaller than $0.2 \lambda$. Fig.~\ref{fig:mlfma_near} demonstrates the controllable accuracy of computing these interactions as a function of $p$, the expansion coefficient for our wideband MLFMA scheme. As is evident from this figure, one can control the accuracy to very fine precision.}


\subsection{Accuracy of CC-CFIER}
In the first set of numerical results, we aim to compare the accuracy and high order nature of the proposed approaches for the analysis of EM scattering against an analytical solution, as well as the number of iterations required by the GMRES solver to reach the prescribed tolerance. To this end, we consider a sphere of diameter 8$\lambda$ that is modeled using an initial control mesh comprising of 642 vertices and 1280 faces. We consider two meshes generated by refining the initial control mesh once and thrice, respectively, using Loop subdivision. \emph{Note}, unlike typical mesh refinement, under the rules of subdivision, the limit surface that \emph{all} meshes point to is identical. More to the point, all the required numerics are \emph{carried out on the limit surface, NOT the Lagrangian geometric approximation}. This refinement process leads to a coarser sphere of 2,562 vertices and 5,120 faces and a finer one composed of 40,482 vertices and 80,960 faces. The main benefits in refining a mesh is better approximation of the physics on the limit surface. 

\begin{figure}[!htpb]
  \centering
  \resizebox{\columnwidth}{!}{

\includegraphics[width=0.5\textwidth]{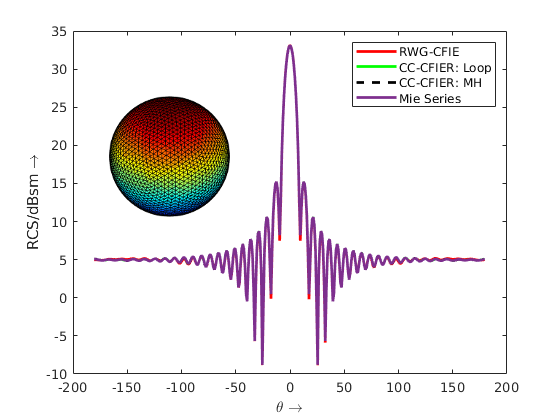}
}
  \caption{Radar cross section of the sphere ($\phi$ = 0 cut).
  \label{fig:rcs_sphere}}  
\end{figure}

In the experiments discussed next, the finer discretization was used with RWG basis (together with a Lagrangian geometry description). We ensured that the surface areas of the Lagrangian mesh agree within 99\% to the subdivision mesh. In Fig.~\ref{fig:rcs_sphere}, we compare RCS data on an 8$\lambda$ sphere \reviewerone{for all three candidates.}

\reviewerone{For the CC-CFIER: MH, we use 1000 MHs leading to 2000 DoF which converges in 7 iterations for a total solve time of 33 seconds; RWG-CFIE requires 122,880 DoF, converges in 36 iterations in 166 seconds, and CC-CFIER: Loop contains 5124 DoF, converges in 7 iterations in 35 seconds.} As is evident from Fig. \ref{fig:rcs_sphere}, the agreement between the three sets of numerical data to analytical solutions is excellent.  In addition, we have analyzed a series of electrically larger spheres. These geometries are obtained via refinement of the initial mesh, such that at any frequency, the edge length is approximately $0.3 \lambda$. The details of these experiments are presented in  Table.~\ref{table2}. As is evident from this table, there is excellent agreement between the proposed methods and analytic data. The convergence of Loop and MH implementations of CC-CFIER is approximately the same as is the total solve time. The approximately four fold compression is not sufficient to affect the overall solve time due to the well-conditioned gram matrix for the sphere.

\begin{table}[!h]
\begin{center}

\resizebox{\columnwidth}{!}{
\begin{tabular}{ c c c c c c}
\hline
 \multirow{3}{*}{Size} &
 \multirow{3}{*}{$N_{L}$/$N_{MH}$} & \multicolumn{2}{c}{CC-CFIER-Loop} & \multicolumn{2}{c}{CC-CFIER-MH} \\
\cmidrule(lr){3-4} \cmidrule(lr){5-6}
& & It./Total Time & $\epsilon^{\infty}$ & It./Total Time & $\epsilon^{\infty}$ \\
\hline
$8\lambda$ & 5124/2000 & 7/0m 35s & 5.99E-4 & 7/0m 33s &  6.32E-4\\

$16\lambda$ & 20484/6000 & 8/4m 31s & 5.99E-4& 8/4m 26s & 9.29E-4 \\

$32\lambda$ & 81924/24000 & 9/25m 42s
 & 2.26E-4 & 9/25m 47s & 2.33E-3\\

\hline
\end{tabular}
}
\end{center}
\caption{Convergence data for a spheres of different diameters: 8$\lambda$-32$\lambda$}
\label{table2}
\end{table}

\begin{table}[h!]
\begin{center}
\resizebox{\columnwidth}{!}{
\begin{tabular}{ c c c c c}
\hline
 \multirow{3}{*}{Size} &
 \multirow{3}{*}{$N_{L}$/$N_{MH}$} & CC-CFIER-Loop & \multicolumn{2}{c}{CC-CFIER-MH} \\
\cmidrule(lr){3-3} \cmidrule(lr){4-5}
& & It./Total Time  & It./Total Time & $\epsilon^{\infty}$ \\
\hline
$8\lambda$  & 5124/2400  & 11/1m 4s & 11/1m 2s & 2.13E-3 \\

$16\lambda$ & 20484/7200 & 12/7m 0s & 12/6m 37s & 5.13E-3 \\

$32\lambda$ & 81924/28000 & 12/38m 29s& 13/36m 10s & 4.44E-3  \\

\hline
\end{tabular}
}
\end{center}
\caption{Convergence data for a bumpy cube of sizes varying from $8\lambda-32\lambda$}
\label{table3}
\end{table}

\subsection{EM Scattering from Complex Objects}
In this section, we provide several examples to demonstrate the viability of using the formulations presented here for EM scattering on complex objects. We do so by comparing our results obtained from CC-CFIER: MH against those obtained using the CC-CFIER: Loop and RWG-CFIE. 

\begin{figure}[!h]
  \centering
  \resizebox{\columnwidth}{!}{
\includegraphics[width=0.5\textwidth]{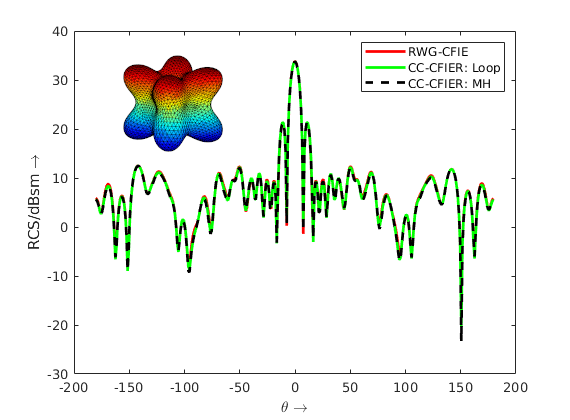}
}
  \caption{Radar cross section of the bumpy cube ($\phi$ = 0 cut).
  \label{fig:rcs_bumpy_cube}}  
\end{figure}

\reviewerone{First, we consider the bumpy cube shown in  Fig.~\ref{fig:rcs_bumpy_cube}, that fits in a 8$\lambda$ $\times$  $8\lambda$ $\times$ $8 \lambda$ box. The number of DoFs for the RWG-CFIE is 122880, converges in 47 iterations for a total of 277 seconds. Whereas, CC-CFIER: Loop and CC-CFIER: MH require 5124 and 2400 degrees of freedom, respectively. Both converge in 11 iterations for a solve time of 60 seconds.} Fig.~\ref{fig:rcs_bumpy_cube} illustrates excellent agreement between all three.

As before, we use mesh refinement to generate electrically larger structures. The results of these runs are presented in Table.~\ref{table3}, specifically, iteration count for CC-CFIER: Loop and CC-CFIER: MH formulation. We report that the iteration count is low, approximately the same for both Loop and MH, and both took approximately the same time for the matrix solve. The  agreement between Loop and the compressed MH system is also excellent. 

\begin{figure}[!h]
  \centering
\includegraphics[width=0.5\textwidth]{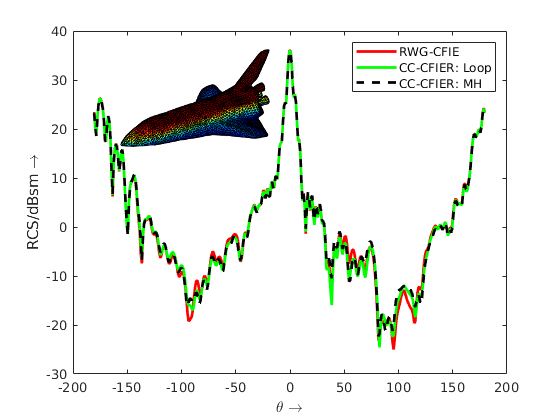}
  \caption{Radar cross section of the shuttle ($\phi$ = 0 cut).
  \label{fig:rcs_shuttle}}  
\end{figure}

Next, we consider a shuttle that that fits in a 20$\lambda$ $\times$ $12.22\lambda$ $\times$ $7.22 \lambda$ box. \reviewerone{The number of DoFs for the RWG-CFIE is 190080, converges in 273 iterations for a total time of 1202 seconds. The CC-CFIER: Loop uses 31684 DoFs, converges in 78 iterations that takes 684 seconds. Lastly, for CC-CFIER: MH uses 6000 DoFs, converges in 39 iterations that takes 311 seconds.} From Fig.~\ref{fig:rcs_shuttle} shows excellent agreement between all three. Again, we refine the geometry to consider electrically larger scatterers, in this case up to 80$\lambda$. Table.~\ref{table4} reports the iteration count, for CC-CFIER: Loop and CC-CFIER: MH basis, as we increase the frequency. We find that the iteration count is stable for both formulation, and they are in excellent agreement. Further, we note the significant compression achieved via MHBs.  
\begin{table}[h!]
\begin{center}
\resizebox{\columnwidth}{!}{
\begin{tabular}{ c c c c c}
\hline
 \multirow{3}{*}{Size} &
 \multirow{3}{*}{$N_{L}$/$N_{MH}$} & CC-CFIER-Loop & \multicolumn{2}{c}{CC-CFIER-MH} \\
\cmidrule(lr){3-3} \cmidrule(lr){4-5}
& & It./Total Time  & It./Total Time & $\epsilon^{\infty}$ \\
\hline
$20\lambda$  & 7942/4000 & 78/11m 24s & 39/5m 15s &2.43E-3  \\

$40\lambda$ & 31684/12000 & 38/30m 57s  & 29/19m 46s  &2.00E-3 \\

$80\lambda$ & 126724/36000 &28/187m 30s & 29/105m 47s & 2.70E-3 \\

\hline
\end{tabular}
}
\end{center}
\caption{Data for shuttle geometries from $20\lambda-80\lambda$.}
\label{table4}
\end{table}

\begin{figure}[!h]
  \centering
  \resizebox{\columnwidth}{!}{
\includegraphics[width=0.5\textwidth]{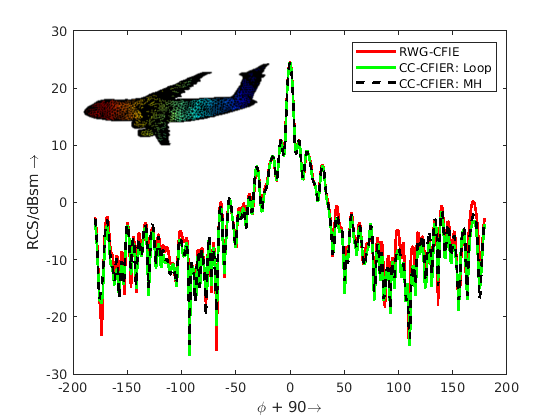}
}
  \caption{Radar cross section of the jet airliner ($\theta$ = 90 cut).
  \label{fig:rcs_jet_airliner}}  
\end{figure}

Finally, we consider a Jet airliner that fits in a 18$\lambda$ $\times$ $17 \lambda$ $\times$ $5 \lambda$ box. In this example, the plane wave propagating in the $\hat{y}$ direction (incident on the nose) and polarized along $\hat{x}$ direction.\reviewerone{ The number of DoFs for the RWG-CFIE is 72768, converging in 243 iterations in 1 hour,  whereas for the CC-CFIER:Loop is 12132 68 iterations in 12 minutes and 45 seconds and the CC-CFIER: MH is 5000 and reaches tolerance within 45 iterations in 5 minutes.} It is evident from Fig.~\ref{fig:rcs_jet_airliner} that all three data sets agree well with each other. In Table.~\ref{table5}, we report the iteration count, for CC-CFIER: Loop and CC-CFIER: MH basis, as we increase the electrical size of the object. We find that the iteration count is stable for both formulation, as well as excellent agreement. Also, note the excellent compression produced by MHBs.

\begin{table}[h!]
\begin{center}
\resizebox{\columnwidth}{!}{
\begin{tabular}{ c c c c c}
\hline
 \multirow{3}{*}{Size} &
 \multirow{3}{*}{$N_{L}$/$N_{MH}$} & CC-CFIER-Loop & \multicolumn{2}{c}{CC-CFIER-MH} \\
\cmidrule(lr){3-3} \cmidrule(lr){4-5}
& & It./Total Time  & It./Total Time & $\epsilon^{\infty}$ \\
\hline
$30\lambda$  & 12132/7000 & 57/14m 59s & 39/9m 45s & 3.97E-3 \\

$60\lambda$ & 48516/21000 & 42/82m 12s  & 39/52m 10s   & 4.78E-3   \\

$120\lambda$ & 194052/63000 & 41/739m 46s & 40/275m 28s &1.34E-2  \\

\hline
\end{tabular}
}
\end{center}
\caption{Data for jetliner geometries from $30\lambda-120\lambda$.}
\label{table5}
\end{table}

\section{Summary}

In this paper, we have presented isogeometric analysis method for subdivision surface; in presenting this approach, we assumed a simply connected structure, used a complete surface Helmholtz decomposition to effect a Calder\'{o}n operator. The latter is possible because the surface is $C^2$ almost everywhere. To evaluate inner-products, which are the main bottleneck for higher order basis functions on higher order surfaces, we use wideband MLFMA to evaluate \emph{all} interactions. Finally, we introduce the notion of manifold harmonics as a means to represent the currents on the surface. These geometry basis can be used for compression of both the manifold and physics on the manifold. We present numerous results using both the subdivision and MH basis, on a collection of electrically large geometries. Two salient points that are evident, (a) subdivision basis are excellent candidates for analysis and (b) MHB provide a mapping on to the eigen-structures of Debye-potentials on the surface. While one can get the compression expected due to a global eigenstructure, a problem that we have not addressed in this paper is the cost of effecting this transformation \cite{nasikun2021hierarchical}. One avenue in particular that aims to mitigate the costs of the MHT is the use of a set of MHs generated by a point-wise product of a small subset of the original MHB \cite{tensor_MH}. \reviewerone{Furthermore, there are a number of capabilities that are still missing; the two most significant are (a) extension to multiply connected structures due to using an exact Helmholtz decomposition and (b) open structures. We have made significant progress on the former and a paper detailing these has been submitted. The second is only a consequence of the subdivision scheme that we have chosen \cite{IGA_JCP}; indeed, one could enrich the current scheme with RWG like basis sets in a manner similar to \cite{JieDaultShankerChapterSubd} or use a different subdivision scheme that permits edges; this and other features of this method will be addressed in subsequent papers.} 
 
\section*{Acknowledgments}
The  authors  acknowledge  computing  support  from  the  HPC  Center  at
Michigan  State  University,  financial  support  from  NSF  via  CMMI-1725278 and  US  Air  Force  Research  Laboratory  under  contracts FA8650-19-F-1747 and FA8650-20-C-1132.

\bibliographystyle{IEEEtran}
\bibliography{IEEEabrv,Bibliography}

\end{document}